\documentclass[graybox]{svmult}

\RequirePackage{type1cm}        \RequirePackage{makeidx}         \RequirePackage{graphicx}        \RequirePackage{multicol}        \RequirePackage[bottom]{footmisc}\RequirePackage{amssymb}
\RequirePackage{url}

\DeclareFontFamily{U}{mathb}{\hyphenchar\font45}
\DeclareFontShape{U}{mathb}{m}{n}{ <-6> matha5 <6-7> matha6 <7-8>
mathb7 <8-9> mathb8 <9-10> mathb9 <10-12> mathb10 <12-> mathb12 }{}
\DeclareSymbolFont{mathb}{U}{mathb}{m}{n}

\DeclareMathAccent{\abxring}{0}{mathb}{"38}

\DeclareFontFamily{U}{mathb}{\hyphenchar\font45}
\DeclareFontShape{U}{mathb}{m}{n}{ <-6> matha5 <6-7> matha6 <7-8>
mathb7 <8-9> mathb8 <9-10> mathb9 <10-12> mathb10 <12-> mathb12 }{}
\DeclareSymbolFont{mathb}{U}{mathb}{m}{n}

\RequirePackage{newtxtext}       \usepackage{enumitem}

\RequirePackage{calrsfs}

\RequirePackage{amsfonts}

\RequirePackage[most]{tcolorbox}

\RequirePackage{bbm}

\RequirePackage{graphicx}

\RequirePackage{booktabs}

\RequirePackage{url}

\RequirePackage{hyperref}

\RequirePackage{color}

\RequirePackage{wasysym}
\RequirePackage{listings}
\RequirePackage{caption}
\RequirePackage{subcaption}

\RequirePackage{stmaryrd}

\RequirePackage{lineno}

\RequirePackage{xparse}

\RequirePackage{float}

\RequirePackage{multirow}

\RequirePackage{cellspace}

\RequirePackage{shellesc}

\RequirePackage[ruled,vlined,linesnumbered]{algorithm2e}

\SetCommentSty{xCommentSty}

\RequirePackage{tikz}
\usetikzlibrary{calc, automata, positioning, arrows, external, matrix, shapes.misc}
\tikzset{cross/.style={cross out,thick,draw=black,minimum size=2*(#1-\pgflinewidth), inner sep=0pt, outer sep=0pt},
cross/.default={2.5pt}}

\RequirePackage{pgfplots}
\pgfplotsset{compat=1.16}
\usepgfplotslibrary{fillbetween}
\usepgfplotslibrary{groupplots}
\usepgfplotslibrary{external}

\DeclareMathAlphabet{\pazocal}{OMS}{zplm}{m}{n}
\DeclareMathAlphabet\bpazocal{OMS}{cmsy}{b}{n}

\providecommand{\bR}{\mathbb{R}}

\providecommand{\vertiii}[1]{{\left\vert\kern-0.15ex\left\vert\kern-0.15ex\left\vert #1
    \right\vert\kern-0.15ex\right\vert\kern-0.15ex\right\vert}}

\NewDocumentCommand{\curlii}{sO{}m}
{
	\IfBooleanTF{#1}
    {\dgalext{#3}}
    {\dgalx[#2]{#3}}
}

\NewDocumentCommand{\dgalext}{m}{  \sbox0{    \mathsurround=0pt     $\left\{\vphantom{#1}\right.\kern-\nulldelimiterspace$  }  \sbox2{\{}  \ifdim\ht0=\ht2
    \{\kern-.625\wd2 \{#1\}\kern-.625\wd2 \}  \else
    \left\{\kern-.7\wd0\left\{#1\right\}\kern-.7\wd0\right\}  \fi
}

\NewDocumentCommand{\dgalx}{om}{  \sbox0{\mathsurround=0pt$#1\{$}  \sbox2{\{}  \ifdim\ht0=\ht2
    \{\kern-.625\wd2 \{#2\}\kern-.625\wd2 \}  \else
    \mathopen{#1\{\kern-.7\wd0 #1\{}
    #2
    \mathclose{#1\}\kern-.7\wd0 #1\}}
  \fi
}

\tikzset{
  partial ellipse/.style args={#1:#2:#3}{
    insert path={+ (#1:#3) arc (#1:#2:#3)}
  }
}

\providecommand{\pF}{\pazocal{F}}

\providecommand{\pP}{\pazocal{P}}

\providecommand{\pT}{\pazocal{T}}

\providecommand{\pV}{\pazocal{V}}

\providecommand{\to}{\widetilde{o}}

\definecolor{blackmy}{RGB}{38, 70, 83}

\definecolor{bluemy}{RGB}{39, 125, 161}
\definecolor{greenmy}{RGB}{42, 167, 143}
\definecolor{yellowmy}{RGB}{233, 196, 106}
\definecolor{brownmy}{RGB}{244, 162, 97}
\definecolor{redmy}{RGB}{249, 65, 68}

 \usepackage{xfrac}
\usepackage{color, colortbl}

\usepackage{lineno}

\makeatletter
\newcommand{\thickhline}{\noalign {\ifnum 0=`}\fi \hrule height 1pt
    \futurelet \reserved@a \@xhline
}
\newcolumntype{"}{@{\hskip\tabcolsep\vrule width 1pt\hskip\tabcolsep}}
\def\blfootnote{\gdef\@thefnmark{}\@footnotetext}
\makeatother

\makeindex

\begin{document}

\title*{Nonlinear Schwarz preconditioning for nonlinear optimization problems with bound constraints}
\titlerunning{Nonlinear preconditioning for bound-constrained optimization problems}
\author{Hardik Kothari, Alena Kopani\v{c}\'{a}kov\'{a} and Rolf Krause}
\institute{Hardik Kothari  \at Universit\`{a} della Svizzera italiana, Switzerland, \email{hardik.kothari@usi.ch}
\and Alena Kopani\v{c}\'{a}kov\'{a} \at Brown University, USA,  \email{alena.kopanicakova@brown.edu}\\ Universit\`{a} della Svizzera italiana, Switzerland, \email{alena.kopanicakova@usi.ch}
\and Rolf Krause \at Universit\`{a} della Svizzera italiana, Switzerland, \email{rolf.krause@usi.ch}
}
\maketitle
\abstract*{We propose a nonlinear additive Schwarz method for solving nonlinear optimization problems with bound constraints.
Our method is used as a ``right-preconditioner" for solving the first-order optimality system arising within the sequential quadratic programming (SQP) framework using Newton's method.
The algorithmic scalability of this preconditioner is enhanced by incorporating a solution-dependent coarse space, which takes into account the restricted constraints from the fine level.
By means of numerical examples, we demonstrate that the proposed preconditioned Newton methods outperform standard active-set methods considered in the literature.
 }
\abstract{We propose a nonlinear additive Schwarz method for solving nonlinear optimization problems with bound constraints.
Our method is used as a ``right-preconditioner" for solving the first-order optimality system arising within the sequential quadratic programming (SQP) framework using Newton's method.
The algorithmic scalability of this preconditioner is enhanced by incorporating a solution-dependent coarse space, which takes into account the restricted constraints from the fine level.
By means of numerical examples, we demonstrate that the proposed preconditioned Newton methods outperform standard active-set methods considered in the literature.
 }

\section{Introduction}
\label{sec:intro}
We consider a Lipschitz domain $\Omega \subset \bR^d$, $d=2,3$, and a triangulation $\pT$ on $\Omega$.
Now, we define $\pV = \text{span}\{\phi_p\}_{p\in \mathcal{N}}$ as a Finite Element (FE) space, where $\mathcal{N}$ denotes a set of nodes of the mesh $\pT$.
Furthermore, we introduce the feasible set ${\pF = \{v \in \pV\ |\ \underline{\psi} \leqslant v \leqslant \overline{\psi} \}}$, where $\underline{\psi},\overline{\psi}$ denote the component-wise lower bound and upper bound, respectively.
We consider the following abstract nonlinear minimization problem:
\begin{equation}
  \begin{aligned}
    \text{Find } & v^{\ast} = {\text{arg}\,\min}_{v\in \pF} f(v),
  \end{aligned}
  \label{eq:problem}
\end{equation}
where $f:\pV \to \bR$ denotes a bounded, twice-Lipschitz-continuously-differentiable objective function.
Problems of this type arise in numerous applications, such as contact mechanics~\cite{kothari2022generalized}, or fracture mechanics~\cite{kopanivcakova2020recursive, kopanivcakova2022nonlinear}.

Under certain assumptions on the function $f$, the minimization problem \eqref{eq:problem} can be equivalently rewritten as a nonlinear complementarity problem (NCP).
As the first-order optimality conditions for \eqref{eq:problem} are given as:
\begin{equation}
  \label{eq:NCP}
  \text{Find } v^\ast \in \pV \text{ such that }
  \begin{cases}
    \nabla f(v) \geqslant 0, \qquad \overline{\psi} - v \geqslant 0,\qquad v - \underline{\psi} \geqslant 0, \\
    \langle \nabla f(v), \overline{\psi} - v \rangle = 0, \quad \langle \nabla f(v), v - \underline{\psi} \rangle = 0.
  \end{cases}
\end{equation}
Standard approaches for solving such minimization problems include penalty/augmented Lagrangian methods, interior-point methods, or active-set methods; see~\cite{jorgenocedal2000-04-27} for a detailed overview.
In this work, we focus our attention on Newton-based active-set methods, namely the semismooth Newton method, and the sequential quadratic programming (SQP) Newton method.
Although the active-set methods are fairly efficient, their convergence tends to deteriorate due to three main factors: inability to detect an active-set sufficiently fast; strong and highly unbalanced nonlinearities; and ill-conditioning of the problem.

In the context of unconstrained nonlinear problems, nonlinear additive Schwarz preconditioners have been demonstrated to accelerate the convergence of the Newton methods; see for example Additive Schwarz Preconditioned Inexact Newton (ASPIN)~\cite{cai_nonlinearly_2002}, Restricted Additive Schwarz Preconditioned Exact Newton (RASPEN)~\cite{dolean_nonlinear_2016}, or Nonlinear Elimination Preconditioned Inexact Newton (NEPIN)~\cite{cai2011inexact} methods.
In this work, we aim to extend a class of Schwarz preconditioned Newton methods to solve constrained nonlinear optimization problems.
To the best of our knowledge, there have been only a few attempts to employ Schwarz methods to solve variational inequalities, for instance~\cite{badea2000additive, badea2003convergence, gross2009additive,badea2012one}.
In this work, we introduce a two-level nonlinear additive Schwarz preconditioner for the Newton-SQP method, which ensures that the subdomain and coarse-level corrections remain in the feasible set.

\section{Nonlinear Preconditioning} \label{sec:NRASPN}
We define a residual function $F:\pV \to \pV'$ as the gradient of the original objective function, i.e.,~$F(\cdot) \equiv \nabla f(\cdot)$.
Let $G$ be a nonlinear preconditioner of the residual function $F$, such that in some sense $G$ denotes an approximate inverse of the nonlinear operator $F$, i.e.,~$G\approx F^{-1}$.
This preconditioner is used to define a nonlinearly-preconditioned variational inequality
problem as follows: Find $v^{\ast} \in \pV$ such that
\begin{equation}
  F(G(v)) =: \mathcal{F}(v) \geqslant 0, \quad \langle \mathcal{F}(v), \overline{\psi} - v \rangle = 0, \quad  \langle \mathcal{F}(v), v - \underline{\psi} \rangle = 0,
  \label{eq:prec_res}
\end{equation}
where the solution of~\eqref{eq:prec_res} is the same as the solution of \eqref{eq:problem} and $\eqref{eq:NCP}$.
Please note that the operator $G$ is used as a ``right-preconditioner", since this type of preconditioning does not change the original nonlinear system, and it also avoids the transformation of bound constraints into general inequality constraints.
Generally, the preconditioner $G$ can be constructed implicitly as a fixed-point iteration, i.e., $v^{(k+1)} = G(v^{(k)})$.
In this work, we construct $G$ using a variant of the nonlinear Restricted Additive Schwarz (NRAS) method, termed as NRAS-B method.
Contrary to the standard NRAS method, the NRAS-B method ensures that the bound constraints are not violated by taking the preconditioning step.
Thus, the preconditioner $G$ produces an iterate that remains in the feasible set, i.e., $v^{(k)} \in \pF$, for all $k=1, 2, \ldots$

\vspace{0.3cm}
\noindent\textbf{NRAS-B method:}
We consider a decomposition of the FE space $\pV$ into $n$ overlapping and non-overlapping subspaces, denoted by $\{\pV_i\}_{i=1}^n$ and $\{\widetilde{\pV}_i\}_{i=1}^n$, respectively.
The overlap between the subspaces is controlled by the variable $\delta$, defined as a multiple of the mesh-width $h$ of the underlying mesh $\pT$.
On these subspaces, we define the standard restriction operator $R_i:\pV\to \pV_i$, and the prolongation operator $P_i:\pV_i \to \pV$, where $R_i^{\top} = P_i$.
Similarly, we define a restricted prolongation operator $\widetilde{P}_i:\widetilde{\pV}_i \to \pV$ such that $\sum_{i=1}^n \widetilde{P}_i R_i = I$.

Utilizing the aforementioned decomposition and the transfer operators, we now define the constrained nonlinear optimization problem on each subspace as follows.
For a given initial guess $v_i \in \pF_i$, where $ v_i \mapsfrom  R_i v^{(k)}$:
\begin{equation}
  \text{Find } v_i^{\ast}  =  {\text{arg}\,\min}_{v_i \in \pF_i} f_i(v_i), \label{eq:subproblem}
\end{equation}
where $f_i:\pV_i \to \bR$ denotes a restriction of the function $f$ to the subspace $\pV_i$.
The feasible set associated with the subspace $\pV_i$ is given by ${\pF_i=\{ v_i \in \pV_i \ | \ \underline{\psi}_i \leqslant v_i \leqslant \overline{\psi}_i\}}$.
Here, we point out that the local minimization problems \eqref{eq:subproblem} are solved on overlapping subspaces.
However, the global iterate $v^{(k)}$ is updated using the corrections associated with the non-overlapping subspaces, as in
\begin{equation}
  v^{(k+1)} = v^{(k)} + \alpha \sum_{i=1}^n \widetilde{P}_i (v_i^{\ast} - R_i v^{(k)}),
  \label{eq:RAS_ls}
\end{equation}
where~$v_i^{\ast}$ denotes a solution of~\eqref{eq:subproblem} and $\alpha$ is computed using a line-search strategy.
Solving~\eqref{eq:subproblem} and update rule \eqref{eq:RAS_ls} comprise an iteration of the NRAS-B method.

\vspace{0.3cm}
\noindent\textbf{Two-level NRAS-B method: }
The convergence of additive Schwarz methods is known to deteriorate with an increasing number of subdomains.
In order to achieve algorithmic scalability, it is essential to ensure global information transfer through a coarse space.
In the context of constrained minimization problems, constructing a coarse space is not a trivial task, as one has to ensure that the prolongated corrections from the coarse level provide a sufficient decrease in the objective function $f$, and the updated iterate remains in the feasible set. We construct a coarse-level objective function~$f_0:\pV_0 \to \bR$, where $\pV_0$ denotes a coarse space $\pV_0 \subset \pV$, and $\pT_0$ denotes a mesh associated with the FE space $\pV_0$.
The transfer of information between the coarse level and the original problem is ensured by the prolongation operator $P_0:\pV_0 \to \pV$ and the restriction operator $R_0: \pV' \to \pV_0'$, where $R_0 = P_0^{\top}$.
Moreover, we also employ the projection operator $\Pi_0: \pV \to \pV_0$ in order to transfer primal quantities to the coarse level.

In the context of nonlinear multilevel methods, several approaches for constructing the coarse-level feasible set ${\pF_0 = \{v_0 \in \pV_0 \ | \ \underline{\psi}_0 \leqslant v_0 \leqslant\overline{\psi}_0 \}}$ are considered in the literature~\cite{gelman1990multilevel, gratton2008_inf, kopanivcakova2021multilevel}.
Here, we utilize constraint-projection rules from~\cite{gelman1990multilevel} and construct $\underline{\psi}_0, \overline{\psi}_0$ in a component-wise manner as
\begin{equation}
  \begin{aligned}
    (\underline{\psi}_0)_t & = (v_0)_t +  \underset{j \in \mathcal{N} \cap \ (\abxring{\omega_0})_t} {\max} [ (\underline{\psi} - v^{(k)} )_j ], \\
    (\overline{\psi}_0)_t       & = (v_0)_t +  \underset{j \in \mathcal{N} \cap \ (\abxring{\omega_0})_t} {\min} [ (\overline{\psi} - v^{(k)} )_j ],
    \label{eq:hard_constraints_KV}
  \end{aligned}
\end{equation}
where the symbol $(\cdot)_t$ denotes the value of a function associated with the $t$-th node of the mesh. The support of the basis function $(\phi_0)_t$ is denoted by ${(\omega_0)}_t$.
Now, we can define the optimization problem on the coarse level as follows.
For a given initial guess $v_0 \in \pazocal{F}_0$, where $v_0 \mapsfrom \Pi_0 v^{(k)}$:
\begin{equation}
  \text{Find } v_0^{\ast} = {\text{arg}\,\min}_{v_0 \in \pF_0} \hat{f}_0(v_0).
  \label{eq:subproblem_coarse}
\end{equation}
Please note that minimization problem~\eqref{eq:subproblem_coarse} is defined using an augmented coarse-level objective function~$ \hat{f}_0$, defined as
\begin{equation}
  \hat{f}_0(v_0) = f_0(v_0) + \langle R_0 \nabla f(v^{(k)}) - \nabla f_0(\Pi_0 v^{(k)}), v_0 \rangle,
  \label{eq:coarse_fun}
\end{equation}
where $v^{(k)}$ denotes the current iterate on the fine level.
By adding the first-order consistency term to the objective function $f_0$, we ensure that the gradient of the augmented objective function $\hat{f}_0$ at the first iterate is the restricted fine-level gradient.

We follow an inverted V-cycle approach, where a coarse-level update step is followed by a single step of NRAS-B iteration, i.e., the iterate $v^{(k)}$ is updated as follows:
\begin{equation}
  \begin{aligned}
    {v}^{(k+\sfrac{1}{2})} & = v^{(k)} + \hat{\alpha} P_0(v_0^\ast - \Pi_0 v^{(k)}),                                                     \\
    v^{(k+1)}              & = {v}^{(k+\sfrac{1}{2})} + \alpha \sum_{i=1}^n \widetilde{P}_i (v_i^{\ast} - R_i {v}^{(k + \sfrac{1}{2})}).
  \end{aligned}
  \label{eq:TLRAS_ls}
\end{equation}
The symbol $v_0^{\ast} $ in~\eqref{eq:TLRAS_ls} denotes the solution of the coarse-level minimization problem \eqref{eq:subproblem_coarse}, while $v_i^\ast$ is the solution of the subproblem~\eqref{eq:subproblem} associated with the $i$-th subspace.
The step sizes $ \hat{\alpha} $ and $\alpha$ are again obtained using a line-search algorithm.
Combining solutions of \eqref{eq:subproblem_coarse} and \eqref{eq:subproblem} with update rule \eqref{eq:TLRAS_ls}, we can define an iteration of the TL-NRAS-B method.

\vspace{0.3cm}
\noindent\textbf{Nonlinearly-preconditioned Newton SQP method: }
Finally, we provide a brief description of the nonlinearly-preconditioned Newton-SQP method for bound-constrained optimization problems.
As summarized in Alg.~\ref{alg:spin}, the method consists of two main phases.
First, we invoke a step of the NRAS-B/TL-NRAS-B method in order to obtain an updated iterate~$v^{(k)}$.
Later, we construct a quadratic model~$Q$, which is minimized subject to the bound constraints with the aim of obtaining a new search direction $s^{(k)}$.
In contrast to standard preconditioned Newton methods~\cite{cai2009nonlinear, cai2011inexact}, the minimization of the quadratic model is subjected to pointwise constraints, which ensures that the updated iterates remain in the feasible set $\pF$.

\begin{algorithm}[t]
  \caption{RASPN-B method}
  \label{alg:spin}
  \DontPrintSemicolon
\SetKwInOut{Input}{Input}
  \SetKwInOut{Output}{Output}

  \KwData{$f:\pV \to \bR$,  $\underline{\psi} \in \pV$, $\overline{\psi} \in \pV$, $v^{(0)} \in \pF$, $k \mapsfrom 0$}
  \KwResult{$v^{(k)}$}
  \BlankLine
  \While {$\|[\nabla f(v^{(k)}) ]_\pF\| \geqslant \epsilon_{{\text{atol}}} $}
  {
  For given $v^{(k)}$, find $v^{(+)}$ by using a step of NRAS-B or TL-NRAS-B method \\
Assemble gradient and Hessian: $g \mapsfrom \nabla f(v^{(+)})$, $H \mapsfrom \nabla^2 f(v^{(+)})$ \\ Find $s^{(k)}$ by solving the following constrained quadratic optimization problem
  \vspace*{-3mm}
  $$
    \min_{s^{(k)}} \  Q(s^{(k)}) := \sfrac{1}{2} \langle H s^{(k)},s^{(k)} \rangle + \langle g, s^{(k)}\rangle, \ \text{\ s.~t. } \underline{\psi} - v^{(+)} \leqslant s^{(k)} \leqslant \overline{\psi} - v^{(+)}
  $$\\
  \vspace*{-3mm}
  Find $\alpha^{(k)}$ using a line-search algorithm \\ Update the iterate: $v^{(k+1)} \mapsfrom v^{(+)} + \alpha^{(k)} s^{(k)}$,\ \  $k\mapsfrom k+1$  \\ }
\end{algorithm}

\section{Numerical Experiments}
\label{sec:numerics}
In this section, we investigate the performance of nonlinear Schwarz preconditioners using two constrained minimization problems, namely the ignition and the minimal surface problems.
Both numerical examples are defined on a domain~${\Omega:=(0,1)^2}$ with boundary~${\Gamma=\partial \Omega}$, which is decomposed into four parts: ${\Gamma_l = \{0\} \times [0,1]}$, ${\Gamma_r = \{1\} \times [0,1]}$, ${\Gamma_b = [0,1] \times \{0\}}$ and $\Gamma_t = [0,1] \times \{1\}$.
The discretization is performed using a mesh consisting of $120 \times 120$ uniform quadrilaterals which are further decomposed into triangular elements.
In the case of two-level methods, we also employ a coarser mesh with $30 \times 30$ elements in each direction.

\vspace{0.3cm}
\noindent
\textbf{Ignition: }
We minimize a variant of the Bratu problem, given as:
\begin{equation}
  \begin{aligned}
     & \underset{u \in H^1(\Omega)}{\text{min}} \ \  f_I(u) := \frac{1}{2} \int_{\Omega} \| \nabla u \|^2 - (u e^u - e^u) \ d x - \int_{\Omega} f(x) u \ dx, \\
     & \text{subject to} \ \ \underline{\psi}(x) \leqslant u \leqslant \overline{\psi}(x), \quad \text{a.e. in } \Omega, \qquad  u = 0, \quad \text{ on } \Gamma,
  \end{aligned}
  \label{eq:min_combustion}
\end{equation}
where ${f(x) = (9 \pi^2 + e^{(x_1^2-x_1^3)\sin(3\pi x_2)} (x_1^2-x_1^3) + 6 x_1 - 2)\sin(3\pi x_1)}$.
The bounds are given as ${\underline{\psi}(x) = 0.2 - 8(x_1 - \sfrac{7}{16})^2 - 8(x_2 - \sfrac{7}{16})^2}$ and $\overline{\psi}(x) = 0.5$.

\vspace{0.3cm}
\noindent
\textbf{Minimal Surface: }
This experiment aims to find the minimal surface described by a function $u$ by solving the following minimization problem:
\begin{equation}
  \begin{aligned}
    \min_{u \in H^1(\Omega)}  f_M(u) = \int_\Omega \sqrt{(1 + \|\nabla u \|^2)} \ dx, \qquad \qquad \qquad \\
    \text{subject to } \begin{cases}
                         \qquad \underline{\psi}(x) \leqslant u \leqslant \overline{\psi}(x) \ \  \text{a.e. in } \Omega.      \\
                         u = -0.3 \sin(2\pi x_2) \ \text{on } \Gamma_l, \ \ u = 0.3 \sin(2\pi x_2) \ \text{on } \Gamma_r, \\
                         u = -0.3 \sin(2\pi x_1) \ \text{on } \Gamma_b, \ \ u = 0.3 \sin(2\pi x_1) \ \text{on } \Gamma_t,
                       \end{cases}
  \end{aligned}
  \label{eq:min_surf}
\end{equation}
where the lower bound is prescribed as
${\underline{\psi}(x) = 0.25 - 8(x_1 -0.7)^2 - 8(x_2 -0.7)^2}$ and the upper bound is ${\overline{\psi}(x) = 8(x_1 -0.3)^2 - 8(x_2 -0.3)^2 - 0.4}$.

\vspace{0.3cm}
\noindent\textbf{Setup of the solution strategies:}
For all numerical experiments, we prescribe the overlap $\delta = 3$ and obtain the decomposition into the subdomains using the library METIS.
All considered solution methods terminate if $\|[{\nabla f}]_{\pF}\| \leqslant 10^{-8}$, where ${[{\nabla f}]_\pF = \pP_\pF(x - \nabla f(x)) - x}$ denotes the projected gradient.
Here, the symbol $\pP_\pF$ denotes the projection onto the feasible set $\pF$.
Contrary to the traditional nonlinear RAS methods, the subdomain solvers, coarse-level solvers, and constrained quadratic minimization solvers are terminated using a fairly strict termination criterion, i.e.,~they terminate if $\| [{\nabla f}]_\pF \| \leqslant 10^{-11}$.
Moreover, we employ a line-search method with the Armijo condition for computing the step size in all inner and outer solvers.
The local and the coarse-level solvers for the NRAS-B/TL-NRAS-B methods employ the Newton-SQP method.

\begin{figure}[t]
  \centering
  \begin{tikzpicture}[]
    \begin{groupplot}[
        group style={
group size = 2 by 2,
horizontal sep=3pt,
            vertical sep=30pt,
          },
width=0.56\textwidth,
        height=5.cm,
        ymode=log,xlabel={$\#$ its.}, minor tick num=9, grid=major,
        ymin=1e-9, ymax=1e1,
        xmin=0,label style={font=\small}, tick label style={font=\small}, legend style={font=\small},
        title style={font=\small}
      ]
\nextgroupplot[align=left, title={},ylabel={$ \|[\nabla f ]_\pF\|$},ytick={1e1,1e0,1e-1,1e-2,1e-3,1e-4,1e-5,1e-6,1e-7,1e-8,1e-9}]
      \addplot [very thick, blackmy]table[x=IT, y=PRN, col sep=comma]{./results/bratu/RAS_SD_2_OV_3.csv};
      \addplot [very thick, bluemy, dotted]table[x=IT, y=PRN, col sep=comma]{./results/bratu/RAS_SD_4_OV_3.csv};
      \addplot [very thick, greenmy, dashed]table[x=IT, y=PRN, col sep=comma]{./results/bratu/RAS_SD_8_OV_3.csv};
      \addplot [very thick, brownmy, dashdotted]table[x=IT, y=PRN, col sep=comma]{./results/bratu/RAS_SD_16_OV_3.csv};
      \addplot [very thick, redmy, loosely dashed]table[x=IT, y=PRN, col sep=comma]{./results/bratu/RAS_SD_32_OV_3.csv};
      \legend{ $\#$sbd = 2, $\#$sbd = 4, $\#$sbd = 8, $\#$sbd = 16, $\#$sbd = 32}

      \nextgroupplot[align=left, title={},ylabel={},yticklabels={},ytick={1e1,1e0,1e-1,1e-2,1e-3,1e-4,1e-5,1e-6,1e-7,1e-8,1e-9}]
      \addplot [very thick, blackmy]table[x=IT, y=PRN, col sep=comma]{./results/bratu/RASTL_SD_2_OV_3.csv};
\addplot [very thick,  bluemy, dotted]table[x=IT, y=PRN, col sep=comma]{./results/bratu/RASTL_SD_4_OV_3.csv};
\addplot [very thick, greenmy, dashed]table[x=IT, y=PRN, col sep=comma]{./results/bratu/RASTL_SD_8_OV_3.csv};
\addplot [very thick, brownmy, dashdotted]table[x=IT, y=PRN, col sep=comma]{./results/bratu/RASTL_SD_16_OV_3.csv};
\addplot [very thick,  redmy, loosely dashed]table[x=IT, y=PRN, col sep=comma]{./results/bratu/RASTL_SD_32_OV_3.csv};
\legend{ $\#$sbd = 2, $\#$sbd = 4, $\#$sbd = 8, $\#$sbd = 16, $\#$sbd = 32}

      \nextgroupplot[align=left, title={},ylabel={$\|[\nabla f ]_\pF\|$},ytick={1e1,1e0,1e-1,1e-2,1e-3,1e-4,1e-5,1e-6,1e-7,1e-8,1e-9}]
      \addplot [very thick,  blackmy]table[x=IT, y=PRN, col sep=comma]{./results/minsurf/RAS_SD_2_OV_3.csv};
\addplot [very thick, bluemy, dotted]table[x=IT, y=PRN, col sep=comma]{./results/minsurf/RAS_SD_4_OV_3.csv};
\addplot [very thick, greenmy, dashed]table[x=IT, y=PRN, col sep=comma]{./results/minsurf/RAS_SD_8_OV_3.csv};
\addplot [very thick, brownmy, dashdotted]table[x=IT, y=PRN, col sep=comma]{./results/minsurf/RAS_SD_16_OV_3.csv};
\addplot [very thick, redmy, loosely dashed]table[x=IT, y=PRN, col sep=comma]{./results/minsurf/RAS_SD_32_OV_3.csv};
\legend{ $\#$sbd = 2, $\#$sbd = 4, $\#$sbd = 8, $\#$sbd = 16, $\#$sbd = 32}

      \nextgroupplot[align=left, title={},ylabel={},minor tick num=4,yticklabels={},ytick={1e1,1e0,1e-1,1e-2,1e-3,1e-4,1e-5,1e-6,1e-7,1e-8,1e-9}]
      \addplot [very thick, blackmy]table[x=IT, y=PRN, col sep=comma]{./results/minsurf/RASTL_SD_2_OV_3.csv};
\addplot [very thick, bluemy, dotted]table[x=IT, y=PRN, col sep=comma]{./results/minsurf/RASTL_SD_4_OV_3.csv};
\addplot [very thick, greenmy, dashed]table[x=IT, y=PRN, col sep=comma]{./results/minsurf/RASTL_SD_8_OV_3.csv};
\addplot [very thick, brownmy, dashdotted]table[x=IT, y=PRN, col sep=comma]{./results/minsurf/RASTL_SD_16_OV_3.csv};
\addplot [very thick, redmy, loosely dashed]table[x=IT, y=PRN, col sep=comma]{./results/minsurf/RASTL_SD_32_OV_3.csv};
\legend{ $\#$sbd = 2, $\#$sbd = 4, $\#$sbd = 8, $\#$sbd = 16, $\#$sbd = 32}
    \end{groupplot}

  \end{tikzpicture}
  \caption{Convergence history of NRAS-B (Left) and TL-NRAS-B (Right) methods for the ignition problem (Top) and the minimal surface (Bottom) problem.
    The experiments are performed with an increasing number of subdomains (sbd).}
  \label{fig:convergence_RAS}
\end{figure}
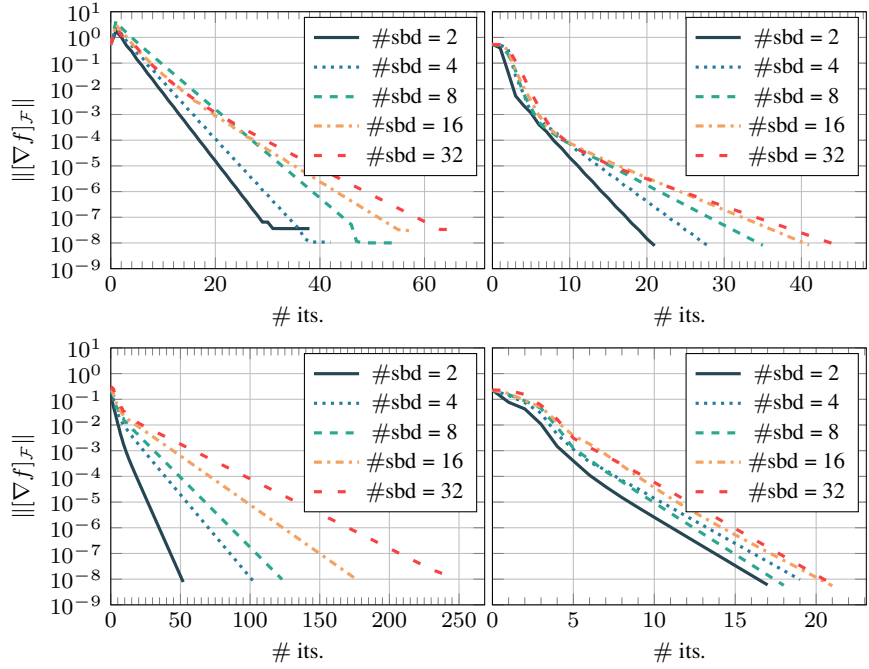

\vspace{0.3cm}
\noindent\textbf{Comparison between NRAS-B and TL-NRAS-B methods: }
The comparison is performed with respect to an increasing number of subdomains.
As we can observe from Fig.~\ref{fig:convergence_RAS}, the standard NRAS-B method requires more iterations to satisfy the termination criterion than the TL-NRAS-B method.
We also notice that the NRAS-B method requires more iterations with an increasing number of subdomains for both problems.
For the TL-NRAS-B, we observe scalable convergence with respect to the number of subdomains for the minimal surface problem.
However, for the ignition problem, the number of iterations grows with the number of subdomains.
This can be attributed to the fact that the coarse grid is not able to represent the constraint function associated with this particular problem sufficiently well.
Hence, the coarse-level nonlinear problems become over-constrained, which amounts to small coarse-grid corrections and insufficient global information transfer.

\vspace{0.3cm}
\noindent\textbf{Comparing RASPN-B method with other methods: }
In this section, we compare the performance of the NRAS-B and TL-NRAS-B preconditioned Newton methods with the semismooth Newton and Newton-SQP methods.
In this study, the NRAS-B and TL-NRAS-B methods employ $16$ subdomains.
We note that the semismooth Newton method linearizes the nonlinearity of the problem and constraints simultaneously, while the Newton-SQP method first linearizes the nonlinearity of the problem and at each Newton iteration a QP problem is solved with constraints.

\begin{figure}[t]
  \centering
  \begin{tikzpicture}[]
    \begin{groupplot}[
        group style={
group size = 2 by 1,
x descriptions at=edge bottom,
horizontal sep=3pt,
          },
        width=0.56\textwidth,
        height=5cm,
        ymode=log,xlabel={$\#$ its.}, minor tick num=4, grid=major,
        ymin=1e-9, ymax=1e0,
        xmin=0,label style={font=\small}, tick label style={font=\small}, legend style={font=\small}
      ]
\nextgroupplot[align=left, title={},ylabel={$\|[\nabla f ]_\pF\|$},ytick={1e0,1e-1,1e-2,1e-3,1e-4,1e-5,1e-6,1e-7,1e-8,1e-9}]
\addplot[blackmy, very thick] table[x=IT, y=PRN, col sep=comma] {./results/bratu/SSN_reduced.csv};
      \label{pgfplots:SSN_reduced}
\addplot[brownmy, dotted, very thick] table[x=IT, y=PRN, col sep=comma] {./results/bratu/NewtonSQP.csv};
      \label{pgfplots:NewonSQP}
\addplot[redmy, dashdotted, very thick] table[x=IT, y=PRN, col sep=comma]{./results/bratu/RASPEN_RAS_SD_16_OV_3.csv};
      \label{pgfplots:RASPEN}
\addplot[greenmy, dashed, very thick] table[x=IT, y=PRN, col sep=comma]{./results/bratu/RASPEN_RASTL_SD_16_OV_3.csv};
      \label{pgfplots:RASTLPEN}

\nextgroupplot[align=left, title={},yticklabels={},ylabel={},ytick={1e-1,1e-2,1e-3,1e-4,1e-5,1e-6,1e-7,1e-8,1e-9}]
\addplot[blackmy, very thick] table[x=IT, y=PRN, col sep=comma] {./results/minsurf/SSN_reduced.csv};
\addplot[brownmy, dotted, very thick]  table[x=IT, y=PRN, col sep=comma] {./results/minsurf/NewtonSQP.csv};
\addplot[redmy, dashdotted, very thick]  table[x=IT, y=PRN, col sep=comma]{./results/minsurf/RASPEN_RAS_SD_16_OV_3.csv};
\addplot[greenmy, dashed, very thick] table[x=IT, y=PRN, col sep=comma]{./results/minsurf/RASPEN_RASTL_SD_16_OV_3.csv};
    \end{groupplot}
    \matrix [ draw, matrix of nodes, anchor = north, node font=\scriptsize,
      column 1/.style={nodes={align=right,text width=0.7cm}},
      column 2/.style={nodes={align=left,text width=1.5cm}},
    ] at ($(group c1r1.south east) + (-2.25, 2.3)$)
    {
      \ref{pgfplots:SSN_reduced} & SS-Newton  \\
      \ref{pgfplots:NewonSQP}    & Newton-SQP \\
      \ref{pgfplots:RASPEN}      & RASPN-B    \\
      \ref{pgfplots:RASTLPEN}    & TL-RASPN-B \\
    };
  \end{tikzpicture}
  \caption{Convergence history of semismooth Newton (SS-Newton),  Newton-SQP, RASPN-B, and TL-RASPN-B methods for the ignition problem (Left) and minimal surface problem (Right).}
  \label{fig:convergence_RASPEN}
\end{figure}
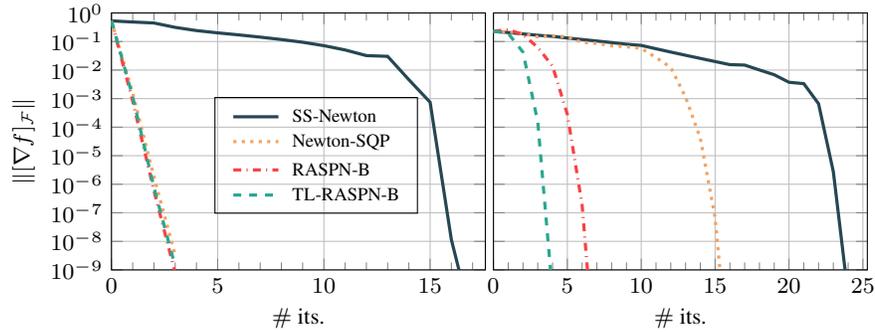

From Fig.~\ref{fig:convergence_RASPEN}, we can see that the Newton-SQP method preconditioned with NRAS-B and TL-NRAS-B method outperforms the semismooth Newton method for both examples.
As the ignition problem is only mildly nonlinear, preconditioned Newton methods as well as the Newton-SQP method manage to satisfy the termination criterion in only $4$ iterations.
In the case of a minimal surface problem, which is more nonlinear, the benefit of preconditioning the Newton method is more evident.
The RASPN-B and TL-RASPN-B methods converge in $7$ and $4$ iterations, respectively.
In comparison, the Newton-SQP and semismooth-Newton methods require $16$ and $24$ iterations to converge, respectively.

\section{Conclusion}
\label{sec:conclusion}

In this work, we presented a nonlinear additive Schwarz preconditioning method for bound-constrained nonlinear optimization problems.
The scalability of the method is enhanced by introducing a coarse level with the first-order consistent objective function and the constraints restricted from the fine level.
The developed Schwarz methods are then employed as a right preconditioner for the Newton-SQP method.
Our numerical results demonstrate that the proposed preconditioners enhance the convergence of the Newton-SQP method and outperform standard active-set Newton methods.
We also show that the two-level preconditioner is algorithmically scalable if a coarse space captures the constraints from the fine level sufficiently well.

\vspace{0.1cm}
\footnotesize{
\noindent\textbf{Acknowledgements }
H.K.~and R.K.~thank the Swiss National Science Foundation (SNSF) and the Deutsche Forschungsgemeinschaft for their support through the project SPP 1962 ``Stress-Based Methods for Variational Inequalities in Solid Mechanics:~Finite Element Discretization and Solution by Hierarchical Optimization" [186407].
A.K.~gratefully acknowledges the support of the SNSF under the project "Multilevel training of DeepONets - multiscale and multiphysics applications" [206745].
Additionally, authors acknowledge the support of Platform for Advanced Scientific Computing through the project ``FraNetG:~Fracture Network Growth" and SNSF through the project ML2 - Multilevel and Domain Decomposition Methods for Machine Learning" [197041].
}
\bibliographystyle{spmpsci}
\bibliography{./extracted.bib}
\end{document}